# Innovative distributed maintenance concept: From the design to cost optimisation

**M. Di Mascolo[1], Z. Simeu-Abazi[2], R. A. Djeunang Mezafack[3]**

[1,2,3] Grenoble Alpes, CNRS, Grenoble INP, G-SCOP, 38000 Grenoble, France

*maria.di-mascolo@grenoble-inp.fr*
*zineb.simeu-abazi@grenoble-inp.fr*
*rony.djeunang-mezafack@grenoble-inp.fr*

## ABSTRACT

This study proposes an integrated heuristic framework for the strategic optimization of distributed maintenance operations in geo-distributed production systems (GDPS). It introduces a dual-entity maintenance structure comprising a Centralized Maintenance Workshop (CMW) and a Mobile Maintenance Workshop (MMW), aimed at minimizing total long-term maintenance costs. The cost function incorporates transport, operations, and downtime penalties, optimized via a two-stage algorithmic approach: a Maintenance Planning Algorithm (MPA) based on predictive maintenance scheduling, and a Long-term Heuristic Scheduling Algorithm (LHSA) addressing a capacitated vehicle routing problem with time windows (CVRPTW). A novel contribution includes a heuristic for CMW location determination using the weighted barycentre of site failure probabilities and a discrete selection of MMW capacities. Mixed Integer Linear Programming (MILP) and a divide-and-conquer heuristic are utilized to handle the NP-hard nature of the problem. Experimental validation using Weibull-distributed failure data and various cost scenarios demonstrates that the proposed Optimised Maintenance and Capacitated Routing (OMCR) framework can reduce lifecycle maintenance costs by up to 50%, with increased scalability for systems exceeding 30 GDPS. The framework is applicable to sectors requiring high availability and centralized servicing, including aerospace, railway, and energy industries.

## 1. INTRODUCTION

Effective maintenance scheduling has always led to a significant improvement in the reliability of industrial systems (Sedghi et al., 2021). It provides the timing of maintenance tasks and the allocation of a set of resources (operators, tools and spare parts). Fortunately, Industry 4.0 technologies (Internet of Things, Artificial Intelligence, Big Data, Digital Twin, etc.) make it possible to anticipate failures in production equipment and offer the possibility of scheduling and managing maintenance operations in an increasingly intelligent manner and in real time (Gopalakrishnan et al., 2022). However, with the rise of global competition in recent decades, manufacturing companies face a highly cost-sensitive market (Saihi et al., 2022). Moreover, maintenance costs can represent between 15% and 70% of total production expenses (Sleptchenko et al., 2019). Therefore, optimising maintenance-related costs is a major issue for companies that want to stay ahead of the competition.

### 1.1. Motivation

Maintenance costs are associated with the resources required for scheduling operations. Obviously, the more equipment there is to maintain, the higher the maintenance costs, especially if the equipment is geographically dispersed. This study focuses on the maintenance of production sites that are located in different places, with equipment that is in use. The simplest approach to organising maintenance for geographically distributed sites is to have each site using its own resources in a decentralised system (Razavi Al-e-hashem et al., 2022). However, such an approach may be more expensive than a centralised system where all sites share the same resources. This paper explores the possibility of centralisation through the concept of distributed maintenance (Manco et al., 2022). The main challenge in this case is resource allocation, as geo-distributed sites share the same maintenance resources (Zhang and Yang, 2021).

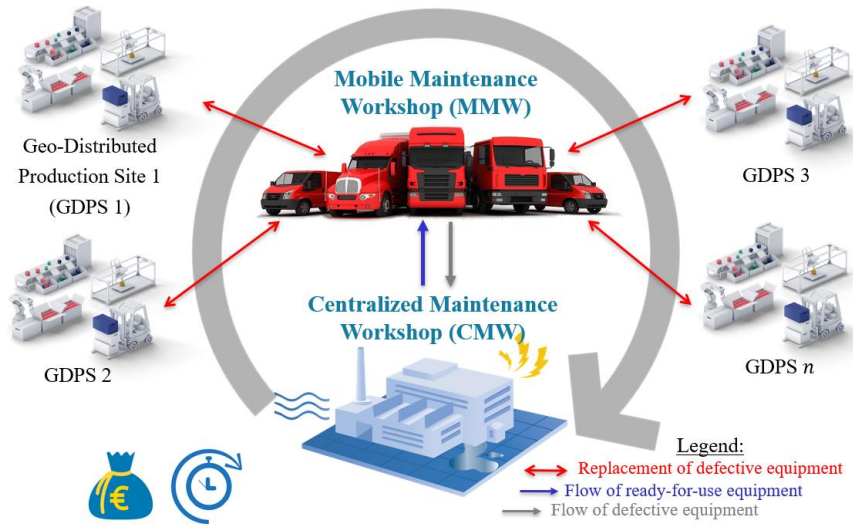

Figure 1: Structure of a Distributed Maintenance, with $n$ GDPS

## 1.2. Problem statement

Distributed maintenance involves a Centralised Maintenance Workshop (CMW) that pools resources, and a Mobile Maintenance Workshop (MMW) that acts as a physical link between Geo-Distributed Production Sites (GDPS), as shown in Figure 1. The MMW follows a predetermined schedule and visits each site to transport spare parts and operators for maintenance operations. The aim of distributed maintenance is to reduce costs by sharing resources. However, centralisation raises the issue of efficient resource allocation. In addition, the combination of scheduling maintenance operations and routing resources to geo-distributed sites is a difficult bi-objective (NP-Hard) optimisation problem.

Industrial applications of distributed maintenance arise in industries where the distance between geo-distributed sites is not too great (Djeunang Mezafack et al., 2022). For example, in the railway sector, several locomotives share the same maintenance workshop for preventive actions (Hani et al., 2007). In the oil & gas industry, oil platforms (onshore or offshore) are geographically distributed according to the sources of raw materials. A centralised platform manages and carries out maintenance operations. Similarly, in the aviation industry, defective aircraft parts are replaced directly on-site without transporting the aircraft. Afterwards, a centralised workshop is needed to diagnose the origin of the failures and repair them (Sanchez et al., 2020). In other applications, a third party maintains distributed facilities owned by different companies.

There are several approaches in the literature that could be useful in the implementation of distributed maintenance. Most of them deal with scheduling maintenance operations without taking mobility into account, or optimising vehicle routes. On the one hand, scheduling is a well-known problem in maintenance management (Valet et al., 2022). The major difficulty is to find the right number of maintenance operations to perform during a time horizon. Too many operations would lead to high equipment idle time and too few operations would increase the probability of equipment failure. On the other hand, routing optimisation is a familiar problem in operations research. It is a combinatorial optimisation problem generally known as VRP (Vehicle Routing Problem) or TSP (Travelling Salesman Problem). The main difficulty encountered is that this is an NP-hard problem, which means that the optimisation time is exponential as a function of the number of sites studied (Konstantakopoulos et al., 2020).

Some papers in the literature allow the implementation of distributed maintenance by combining the two approaches mentioned above (scheduling of operations and vehicle routing). But these classic approaches only optimise costs from an operational point of view (daily). From a strategic point of view (yearly), no study has yet shown the influence of CMW location and MMW capacity on maintenance costs. These parameters are generally considered to be fixed, otherwise the joint optimisation of scheduling and routing would require exponential computing time using standard methods. In addition, instead of a short-term schedule (days), a long-term schedule (years) is needed to obtain a good estimate of maintenance costs during optimisation. Based on this observation, this paper proposes an approach to:

i. perform long-term maintenance scheduling and vehicle routing

ii. design the two main elements of distributed maintenance: CMW location and MMW capacity.

## 1.3. The contributions

The first objective is to provide a novel heuristic allowing to address the computation time problem related to the long-term aspect. Indeed, the majority of current studies that aim to reduce the computation time propose a short cyclic scheduling (Manco et al., 2022). Then, for a larger time horizon, the proposed cycle is repeated. Although this method guarantees high availability of production sites, maintenance costs remain high because the number of operations to be carried out over the time horizon is overestimated. In contrast to known approaches, this paper proposes a heuristic that enables the study of distributed maintenance costs over a wider operating time horizon.

As part of our work on distributed maintenance, we used two metaheuristics for just routing problems, namely: the genetic algorithm and simulated annealing (Ndiaye I. 2014). Given the effectiveness of the latter in combinatorial problems such as the traveling salesman problem or the vehicle routing problem, we applied them to our problem of optimizing transport and logistics costs in the preventive maintenance of a multi-site system. The comparison of the two methods showed that the genetic algorithm generally remains more efficient because it offers a lower maintenance cost. The comparison of computation times also gives the advantage to the genetic algorithm.

The approach used in this paper is based on predictive maintenance which, in combination with MILP (Mixed Integer Linear Programming), identifies for each GDPS the best time to preventively replace equipment and determines the optimal routing of vehicles. Additionally, the second objective is to optimize the location of CMWs and the capacity of MMWs.

Current studies consider that the CMW should be located close to one of the GDPS for optimal cost (Simeu-Abazi and Gascard, 2020). But the more GDPS there are, the more difficult it is to choose the best one. Thus, what is missing is a precise description of the CMW location to be chosen to reduce maintenance costs independently of the number of GDPS involved. In this research, a second heuristic proposes to position the CMW at the weighted barycentre of the GDPS. After determining the location of the CMW, the capacity of the MMW is chosen from a predefined set. Indeed, the number of vehicle types is relatively small and can be grouped into three main categories (light, medium and heavy vehicles).

Following the introduction in Section 1, a review of the literature is presented in Section 2. The next section proposes a problem formulation with relevant assumptions. Section 4 presents the general optimisation framework and the detailed process steps. Then, Section 5 implements the proposed method in a case study through experiments. Section 6 presents and analyses the results obtained. The last section concludes this study and provides some perspectives for future research.

## 2. LITERATURE REVIEW

The simplest approach towards dealing with GDPS maintenance is to wait for failures to occur before carrying out corrective maintenance. In this case, it would be sufficient to have a list of the sites affected by the failures and to find the best path for the workforce to access the sites. This first problem is a so-called joint optimisation of scheduling and workforce routing for the restoration of GDPS (Yulong et al., 2019). (Gupta, 2003) proposed a simulated annealing algorithm, (Drake et al., 2020) a genetic algorithm and (Allaham and Dalalah, 2022) a MILP to maximise the amount of work and minimise the total distance travelled by the workforce. (Cakirgil et al., 2020) combined a MILP with a variable neighbourhood search to complete the highest priority tasks earlier. In order to consider the possible delay of these corrective maintenance interventions, (Hedjazi et al., 2019) developed a multi-agent system maximising the availability of the facilities. However, in most companies, every second that a piece of equipment is down represents a significant loss. Furthermore, corrective maintenance cannot be applied in a long-term strategy, as operations are performed after failures. It would be interesting to be able to act before failures occur.

Preventive maintenance enables proactive action to anticipate failures and improve the availability of the production sites. (Tang et al., 2007) proposed an adaptive memory tabu search, (López-Santana et al., 2016) a MILP, (Fontecha et al., 2020) a LP-based split heuristic, (Nguyen et al., 2019) a combination of a local search genetic algorithm (LSGA) and branch and bound method, and (Jia and Zhang, 2020) a simulated annealing algorithm for routing a set of crews to perform the planned maintenance operations at a near-minimum expected cost per unit time. To consider the possibility of reducing costs by centralising resources, (Simeu-Abazi and Ahmad, 2011) proposed a modular approach based on Petri nets and (Wang and Djurdjanovic, 2018) a discrete event simulation. However, not all of the above approaches consider the limited capacity of transport vehicles. Thus, (Allaham and Dalalah, 2022) proposed a MILP to introduce transport constraints. It is interesting to note that in the case of preventive maintenance, scheduling is cyclical and only allows maintenance costs to be assessed over a short time horizon. Indeed, preventive maintenance cycles do not consider the degradation of equipment over time.

Predictive maintenance aims to use equipment degradation parameters to establish an appropriate maintenance schedule. The particularity of this approach in the literature is the use of prognostic information to maximise the availability of production sites and minimise the total distance travelled to reach these sites (Camci, 2015). (Rashidnejad et al., 2018) proposed a genetic algorithm for scheduling predictive maintenance operations. To consider the transportation capacity constraints, (Si et al., 2022) proposed a MILP. Based on the real-time machine degradation, it is possible to estimate the assets' failure rate and establish a time-varying maintenance cost function to quantify the trade-off between early and delayed maintenance. In addition, (Manco et al., 2022) also proposes to centralise predictive maintenance operations to reduce costs.

A summary of this literature review is presented in Table 1 (see (Djeunang Mezafack 2023) for more details). The preceding analysis shows that predictive maintenance is the most appropriate strategy for long-term scheduling. The constraints considered concern on the one hand the MMW (vehicle capacity and long-term scheduling) and on the other hand the CMW (centralisation of maintenance and choice of the geographical location of the depot). It can be noted that some studies propose exact resolution methods, but they are far from reality and not very useful for the industry. Other studies use heuristics and metaheuristics to deal with more constraints. However, the use of exact methods is still appropriate for simple study cases (generally less than 10 GDPS and one maintenance operation per site to be performed periodically).

The organisation of distributed maintenance is an NP-hard problem, as explained previously (Manco et al., 2022), and the use of heuristic or metaheuristic methods remains the best compromise between calculation time and good approximation. Although existing approaches attempt to optimise maintenance costs, linked to resources, transport and breakdowns, they only solve operational problems, i.e. over a short time horizon (days). What is missing is a more strategic approach (years) to cost optimisation. In addition, current approaches do not consider the influence of CMW location and MMW capacity on maintenance costs. This study therefore aims to fill this gap based on the current predictive maintenance strategy.

## 3. Problem formulation

The objective of this paper is to propose a decision support tool for the choice of the geographical location of a CMW and the capacity of the MMW through long-term scheduling. Let's consider a set of N heterogeneous GDPS. Each GDPS has a piece of equipment that is subject to uncertain failures.

With a schedule, MMW is responsible for transporting maintenance resources (spare parts and tools) to visit all GDPS within a given time horizon. CMW monitors the condition of each piece of equipment and stores spare parts. The MMW starts his route from the CMW, with a limited spare parts capacity, and visits the GDPS following the optimal scheduling. When the MMW reaches a GDPS, the GDPS equipment is replaced with a spare part. The described transport and maintenance network are established to maximise the operational availability of production equipment installed in the GDPS. However, it is essential to ensure the minimisation of transportation and maintenance costs incurred by this network. Therefore, a cost evaluation method should be proposed to ensure that these costs are kept to a minimum.

As developed in Section 2, there are several methods for evaluating costs in a distributed maintenance context. However, the majority of these methods take a short-term approach, assuming that transportation and maintenance are cyclical. In a cyclical framework, it would be sufficient to assess the minimum costs for one cycle and multiply the result by the number of cycles. Unfortunately, this approach can lead to either an overestimation or underestimation of costs, given that transportation and maintenance needs vary from one cycle to another. Hence, it is necessary to propose a method that improves the long-term cost evaluation. Additionally, the choice of the geographical location of the CMW (Central Maintenance Workshop) and the capacity of the MMW requires a study of their impact on long-term costs. In this problem, we have to manage not only the application of maintenance by replacement in each of the production sites (GDPS), the storage capacity of the mobile workshop (MMW) and its routing but also the positioning of the central workshop (CMW). The main hypotheses of the problem are summarized as follows:

### 3.1. Assumptions about GDPS

i. GDPS are N geographically dispersed sites, each with a single type of equipment

ii. All pieces of equipment are mutually independent from the point of view of failures. This means that the state of health of one piece of equipment does not affect that of another. We therefore consider that it is possible to optimise the number of maintenance operations for each site separately.

iii. The probability of a piece of equipment on site i failing at time t is equal to the cumulative distribution function $F_i$ $(T \leq t)$ where T is a random variable of the time to failure.

iv. Each site i is subject to a tight time window $[a_i, b_i]$ outside of which a maintenance operation cannot be carried out.

v. In the event of an unexpected failure, the site waits until one of the vehicles arrives to replace the faulty equipment.

### 3.2. Assumption about CMW and MMW

i. CMW has a spare parts depot with unlimited capacity.

ii. MMW is a fleet of m homogeneous vehicles, each with a limited capacity Q.

Table 1: Classification of pap

| Problem definition | Authors | Correcti |
|---|---|---|
| Joint optimisation of workforce scheduling and routing for restoring critical geo-distributed assets | (Gupta, 2003); (Hedjazi et al., 2019); (Yulong et al., 2019); (Cakirgil et al., 2020); (Drake et al., 2020); (Allaham and Dalalah, 2022) | ● |
| Scheduling technicians for planned maintenance of geographically distributed equipment | (Tang et al., 2007) (López-Santana et al., 2016); (Fontecha et al., 2020); (Nguyen et al., 2019); (Jia and Zhang, 2020) | |
| Joint Optimisation of Preventive Maintenance, Spare Parts Inventory and Transportation Options for Systems of Geographically Distributed Assets | (Simeu-Abazi and Ahmad, 2011); (Wang and Djurdjanovic, 2018) | |
| Multitask scheduling of geographically distributed maintenance tasks | (Allaham and Dalalah, 2022) | |
| Maintenance scheduling of geographically distributed assets with prognostics information | (Camci, 2015) ; (Rashidnejad et al., 2018) | |
| Maintenance grouping and technician routing problem of geographically distributed production systems | (Si et al., 2022) | |
| Maintenance management for geographically distributed assets | (Manco et al., 2022) | |
| **Design and Optimisation of a Distributed Maintenance** | **This paper** | |

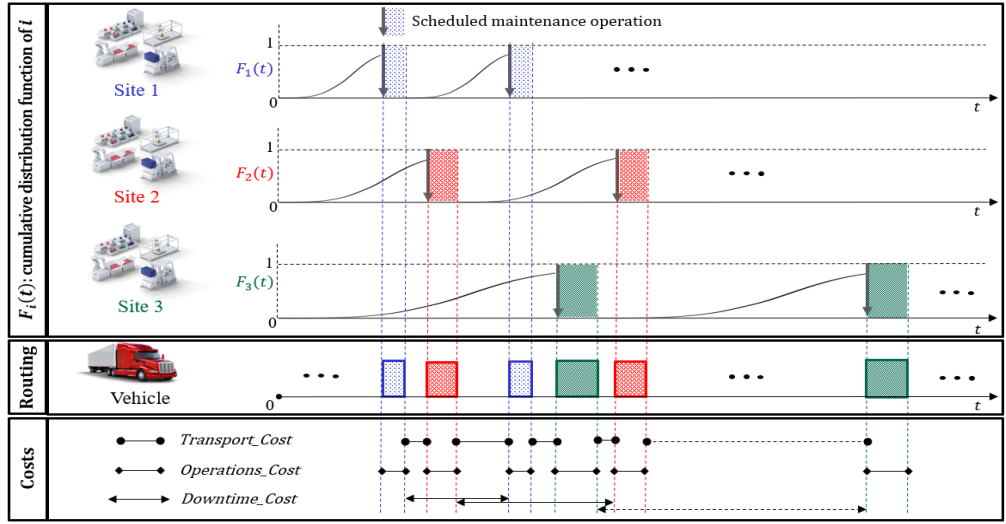

Figure 2: Illustration of distributed maintenance costs for three GDPS

### 3.3. Assumption about vehicle routing

i. A vehicle carries spare parts and one operator, whose role is to replace a piece of equipment with a new spare part.

ii. The travel times $t_{ij}$ and distances $d_{ij}$ between sites i and j are deterministic and do not change over the scheduling horizon τ, i.e. we do not consider any perturbation on the travel times

iii. Each time a vehicle leaves a site i, it has a number $y_i$ of spare parts remaining and as soon as the stock is empty, it returns to the depot for resupply. It may happen that a vehicle returns to the depot without the stock being empty, if this makes the overall routing more efficient.

iv. One of the major problems encountered when implementing distributed maintenance as it has been defined is long-term cost optimisation. The costs have three expected causes, as illustrated in Figure 2: transport; (*Transport_Cost*), operations (*Operations_Cost*) and downtime (*Downtime_Cost*).

- *Transport_Cost* is a linear combination of distances and travel times

- *Operations_Cost* is proportional to the number of maintenance operations.

- *Downtime_Cost* is more complex than the previous two costs but can be determined from the work of (López-Santana et al., 2016). Indeed, for a given piece of equipment, the research shows that the time that elapses between two successive maintenance operations generates a cost linked to the probability that the equipment fails i.e. the cumulative distribution function of the time-to-failure. This is the primary reason why the cumulative distribution function of the time-to-failure is chosen as the indicator of maintenance operation criticality in this study. The higher the probability of equipment failure, the more critical it is considered, and the more priority it will be given in the scheduling of operations.

Other parameters, such as Fussell-Vesely importance (Meng, 2000) or Birnbaum importance (Wu and Coolen, 2013), could have been used as indicators of criticality, but we chose the simplest parameter for this study, as it had also demonstrated its relevance in the context of distributed maintenance in previous studies. Therefore, from equation (31) in (López-Santana et al., 2016), it is sufficient to have the end time of one operation and the start time of the next operation to calculate *Downtime_Cost*. The sum of these three different costs gives the distributed maintenance cost (*Total_Cost*) as expressed in equation (1) to be optimised.

$$Total\_Cost = Transport\_Cost + Operations\_Cost$$
$$+ Downtime\_Cost \qquad (1)$$

The main notations used in this study are summarised in Table 2. In the following section, the details of the proposed algorithms and the general framework to solve the problem are presented.

### 4. DESIGN OPTIMISATION

#### 4.1. Optimised Maintenance and Capacitated Routing

As presented above, the expected costs related to distributed maintenance (*Total_Cost*) are: transport costs (*Transport_Cost*), operating costs (*Operations_Cost*) and downtime costs (*Downtime_Cost*). They can be grouped into two categories based on their origin: transportation (*Transport_Cost*) and maintenance (*Operations_Cost* +

*Downtime_Cost*). Cost optimisation can thus be divided into two steps: maintenance optimisation and transportation optimisation.

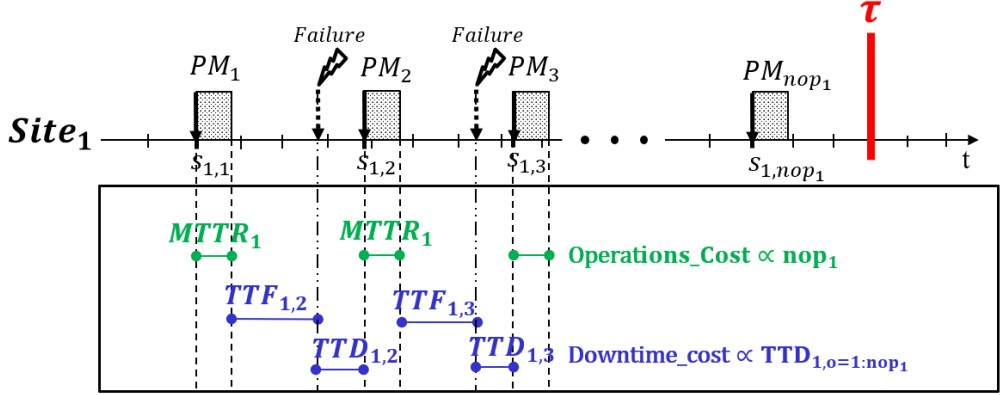

Figure 2: Illustration of distributed maintenance costs for three GDPS

### 4.1.1. Maintenance optimisation

This first step involves minimising *Operations_Cost* and *Downtime_Cost*, and has been extensively studied in scientific literature. Drawing inspiration from the work of López-Santana et al. (2016), the aim is to determine the number of maintenance operations ($nop_i$) to be carried out for each site i within a given horizon $\tau$. For a clearer representation, Figure 3 illustrates the variables that influence these costs for a site.

Let's imagine that for site i, we have planned $nop_i$ predictive maintenance operations. When the equipment at site i is put into operation, it takes a time $TTF_{i,1}$ (Time To Failure) before it fails for the first time. $TTF_{i,1}$ is a random variable with a value ranging from 0 to $s_{i,1}$ (the start time of the first maintenance operation at site i). $TTF_{i,1}$ is associated with the failure probability $F_i$ ($TTF_{i,1} \leq s_{i,1}$). Therefore, site i has a probability $F_i$ ($TTF_{i,1} \leq s_{i,1}$) of being unavailable for a duration of $TTD_{i,1}$ (Time To Downtime). Once it reaches time $s_{i,1}$, whether the equipment is faulty or not, it is replaced/repaired for an average duration of $MTTR_i$ (Main Time To Repair).

The reasoning is similar for subsequent operations until the planned number $nop_i$ of operations is reached. The values of *Operations_Cost* and *Downtime_Cost* can be deduced from Figure 3, considering *N* production sites:

$$Operation\_Cost = \frac{1}{\tau} * \sum_{i=1:N} CR_i . nop_i \qquad (2)$$

$$\begin{aligned} Downtime\_Cost \\ = \frac{1}{\tau} * \sum_{i=1:N} \sum_{o=1:nop_i} CP_i * F_i (TTF_{i,o} \leq s_{i,o}) * TTD_{i,o} \end{aligned} \qquad (3)$$

where $CR_i$ is the cost of a maintenance/replacement operation for the failed equipment on site i and $CP_i$ is the cost of penalty due to the unavailability of site i following a failure.

The decision variables that allow the optimisation of *Operations_Cost* and *Downtime_Cost* through equations (2) and (3) are the elements of the vector [$nop_i$] and the matrix [$s_{i,o}$] given that $\tau$, N, [$CR_i$], [$CP_i$] and [$TTD_{i,o}$] are input parameters.

Regarding the decision variables, we know that, for a site i: $nop_i \in [1, +\infty[$ , $s_{i,o} \in [0, \tau]$, $s_{i,o} < s_{i,o+1}$ and $nop_i$ = number of different values for $s_{i,o}$.

Optimising the above defined costs leads to a nonlinear optimisation problem. Several methods could be used to solve it, such as the Golden-section search, Interpolation methods, Line search, Nelder–Mead method, etc, and they lead to similar results. In this work, we use the Golden-section search method (Chang, 2009). We refer to this first optimisation phase as MPA (Maintenance Planning Algorithm). This algorithm enables to predict the optimal number of maintenance operations $nop_i$ for each site i, and also gives the values of start times $s_{i,o}$ and a set of time windows $[e_{i,o} ; l_{i,o}]$

Table 2: Notation of this study

| Parameters | Description | Unit | Input | Output |
|---|---|---|---|---|
| **GDPS parameters** | | | | |
| $\tau$ | Scheduling horizon | $months$ | • | |
| $N$ | Number of production sites | - | • | |
| $R$ | Distribution radius of the sites | $km$ | • | |
| $t_{ij}, d_{ij}$ | Travel times and distances between sites $i$ and $j$ | $h, km$ | • | |
| $MTTR_i, CR_i$ | Duration and cost of a maintenance/replacement operation | $h, \$$ | • | |
| $CP_i$ | Cost of penalties for unavailability of equipment on site $i$ | $\$/h$ | • | |
| $F_i$ | Cumulative distribution function of the time to failure of equipment on site $i$ | - | • | |
| $nop_i$ | Number of maintenance operations performed on site $i$ | - | | • |
| $n$ | Sum of maintenance operations on all sites | - | | • |
| $s_{i,o}, [e_{i,o}; l_{i,o}]$ | Start time and time window of operation $o$ on site $i$ | $h$ | | • |
| $TTF_{i,o}$ | Time To Failure of the piece of equipment on site $i$ before operation $o$ | $h$ | | • |
| $TTD_{i,o}$ | Time To Downtime of the piece of equipment on site $i$ before operation $o$ | $h$ | | • |
| $A_i$ | Average availability of a production site $i$ during the scheduling horizon | - | | • |
| **CMW parameters** | | | | |
| $\Omega$ | Set of possible CMW positions | - | • | |
| $\Phi_k$ | CMW position at iteration $k$ of the proposed general framework | - | | • |
| **MMW parameters** | | | | |
| $Q$ | Transport capacity of a vehicle | - | • | |
| $CD$ | Unit transport cost per $km$ per unit of capacity | $\$/km$ | • | |
| $CT$ | Unit transport cost per $hour$ | $\$/h$ | • | |
| $m$ | Number of vehicles | - | | • |
| $y_i$ | Number of spare parts remaining in stock after the site visit $i$ | - | | • |
| $x_{ij}$ | Binary decision variable indicating whether a vehicle crosses an arc $(i, j)$ in the optimal scheduling | - | | • |
| **Other notations** | | | | |
| $CMW$ | Centralized maintenance workshop | - | - | - |
| $MMW$ | Mobile maintenance workshop | - | - | - |
| $OMCR$ | Optimized Maintenance and Capacitated Routing | - | - | - |
| $MPA$ | Maintenance planning algorithm | - | - | - |
| $LHSA$ | Long-term Heuristic Scheduling Algorithm | - | - | - |



### 4.1.2. Routing optimisation

The second objective is to minimise Transport_Cost. This problem is quite similar to a Capacitated Vehicle Routing Problem with Time Windows (CVRPTW) in operations research, but it is within the framework of long-term scheduling rather than short-term scheduling. In a classical CVRPTW, a fleet of homogeneous vehicles has to serve customers with known demand and opening hours (Konstantakopoulos et al., 2020). Approaches found in the literature focus on finding the best routes for transport vehicles in such a way that all sites/customers are visited once and only once. In the application of such methods, the optimal route is repeated once all sites have been visited. In the context of distributed maintenance, automatically repeating visit cycles can lead to either an overestimation or underestimation of costs, as mentioned earlier. This is justified by the fact that the maintenance demand for production sites varies over a long-time horizon (several years) according to the failure distribution law. Therefore, one of the innovations of this study will be to propose an algorithm that improves the estimation of maintenance costs for optimal availability of production sites. We will refer to this algorithm as the "Long-term Heuristic Scheduling Algorithm" (LHSA).

In the input of this second algorithm (LHSA), we would have the output of the first algorithm (MPA), i.e., the vector [$nop_i$] and the matrix [$s_{i,o}$]. Additionally, LHSA takes as input a tolerance interval [$e_{i,o}$ ; $l_{i,o}$] for the start time of operations for each site, referred to as "time windows" i.e. $s_{i,o} \in$ [$e_{i,o}$ ; $l_{i,o}$] such that $s_{i,o} - e_{io} = l_{io} - s_{io}$ . The objective is to modify the values of the elements of [$s_{i,o}$] in the time windows [$e_{i,o}$ ; $l_{i,o}$] to achieve a maintenance operation scheduling that guarantees minimal transportation costs.

Let's consider as other input of LHSA a complete directed graph G=(V,A), where V={0,1,2,···,N} is a set of nodes with a depot 0, and $V_s = V\backslash\{0\}$ a subset of nodes.

A={(i, j) : i , j $\in$ V} represents the set of links between all pairs of nodes. The set of vehicles is defined by K={1,2,…,m}, each with a capacity Q. Each site I $\in V_s$ is associated with an on-site service time $MTTR_i$. Non-negative travel times $t_{ij}$ and distances $d_{ij}$ are associated to each arc (i, j) $\in$ A.

Each site I $\in V_s$ has $nop_i$ times a maintenance operation over the scheduling horizon τ. The MPA sub-algorithm provides time windows for each operation [$e_{i,o}$ ; $l_{i,o}$]:o=1,2,…, $nop_i$. An aggregate set of nodes V'={0,1,2,…,n} is therefore considered, where 0 is the depot and n=$\sum nop_i$ represents the sum of all maintenance operations on the horizon τ. Then, an auxiliary directed graph G' = (V', A' ) is defined, where

A'={(i', j' ) : i', j' $\in$ V' } denotes the set of arcs. For each arc (i', j') $\in$ A' the equivalent arc (i, j) $\in$ A can be found such as $t_{i'j'} = t_{ij}$ and $d_{i'j'} = d_{ij}$.

The problem is therefore to determine the optimal routing between the maintenance operations i' such that:

i.   Each operation i' $\in$ V'\{0} is performed exactly once.

ii.  A vehicle cannot transport spare parts over its capacity Q

iii. The time window [$e_{i'}$ ; $l_i$] of operation i' is equivalent to a single time window [$e_{i,o}$ ; $l_{i,o}$] : i $\in V_s$, o=1,2,…,$nop_i$, provided by the MPA sub-algorithm and vice versa.

In the remainder of this paper, the index i is used, instead of i' or o to refer to each maintenance operation. This formulated problem is NP-hard, requiring exponential computation time. Only small instances can be solved analytically. An analytical model is thus defined for solving small instances and a heuristic for larger instances of the problem. A mixed integer linear programming (MILP) model is first chosen, which is the most widely used in the literature (Borcinova, 2017).

Next, a divide-and-conquer algorithm (Mariescu-Istodor et al., 2021) is adapted and implemented to solve the computational time problem, as illustrated in Figure 4.

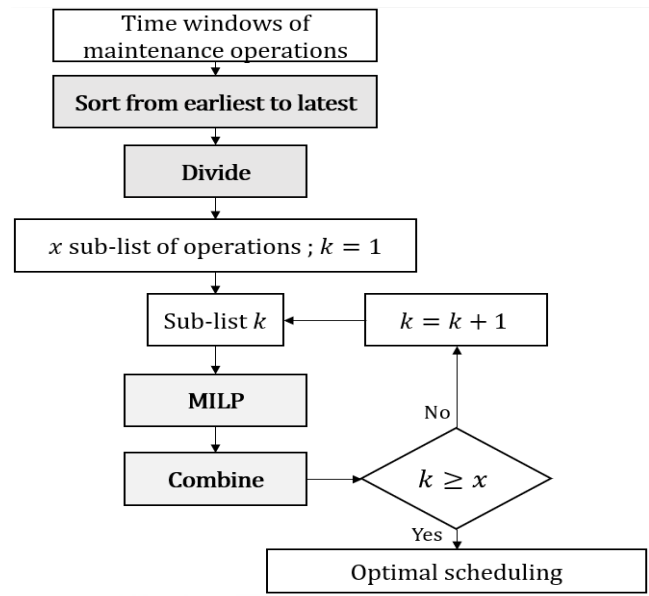

Figure 4 Divide-and-Conquer algorithm, adapted to the LHSA sub-algorithm.

This method should make it possible to deal with the computation time by dividing the list of all maintenance operations into ordered sub-lists. First, the large list of maintenance operations is sorted from oldest to newest and divided into smaller sub-lists. The maximum number of elements in a sub-list is equal to Q. Then each sub-list is solved by the MILP and the results are combined to obtain the maintenance scheduling.





The focus is now on the MILP formulation. A binary decision variable $x_{ij}$ is defined to indicate whether a vehicle crosses an arc (i,j) in the optimal scheduling. A vehicle arrives for an operation i at a time indicated by $s_i$ and with a load $y_i$.

The MILP of the LHSA sub-algorithm can be stated as follows:

Minimize

$$Transport_{Cost} = \frac{1}{\tau} \sum_{i=0}^{n} \sum_{j=0}^{n} \left( Q.CD.d_{ij} + CT.t_{ij} \right).x_{ij} \quad (4)$$

Subject to

$$-Q.m \leq -n; \quad (5)$$

$$\sum_{j=1}^{n} x_{0j} - m = 0; \quad (6)$$

$$\sum_{i=0}^{n} x_{ii} = 0; \quad (7)$$

$$\sum_{i=0, i \neq j}^{n} x_{ij} = 1, \forall j \; \epsilon V'\backslash\{0\}; \quad (8)$$

$$\sum_{j=1, i \neq j}^{n} x_{ij} \leq 1, \forall i \; \epsilon \; V'\backslash\{0\}; \quad (9)$$

$$y_i - y_j + (1 + Q).x_{ij} \leq Q, \forall i, j \epsilon V'\backslash\{0\}, i \neq j; \quad (10)$$

$$s_i - s_j + \left(TR_i + t_{ij} + \tau\right).x_{ij} \leq \tau, \forall i, j \epsilon V'\backslash\{0\}, i \neq j; \quad (11)$$

$$1 \leq y_i \leq Q, \forall i \epsilon V'^{\backslash\{0\}}; \quad (12)$$

$$e_i \leq s_i \leq l_i, \forall i \; \epsilon V'^{\backslash\{0\}}; \quad (13)$$

$$x_{ij} \epsilon \{0,1\}, \forall i, j \; \epsilon \; V'^{\backslash\{0\}}; \quad (14)$$

This MILP formulation minimises the transport costs $Transport\_Cost$ through the objective function (1). Constraint (5) represents the minimum number of vehicles required to service all operations. In the linear program, a vehicle is considered as a route that starts from the central maintenance workshop (CMW), visits a certain number of sites, and returns to the central maintenance workshop (CMW), as illustrated in Figure 5. Constraint (6) requires that exactly m vehicles leave the depot. The classical flow constraints (7), (8) and (9) ensure that each vehicle can leave the depot exactly once, and that each maintenance operation is performed only once. In constraint (10), the capacity of the vehicles is defined such that the difference in load of a vehicle between two successive operations i and j does not exceed the demand of j. Constraint (11) ensures that the time between two successive operations i and j does not exceed $MTTR_i + t_{ij}$. Constraints (12), (13) and (14) restrict the upper and lower bounds of the decision variables.

### 4.1.3. Joint optimisation of maintenance and transport

The optimisation approach proposed for the planning and scheduling of distributed maintenance is an iterative process between the MPA and LHSA algorithms, referred to as OMCR (Optimised Maintenance and Capacitated Routing), as illustrated in Figure 6. Thus, the OMCR algorithm consists in solving two different sub-algorithms iteratively. At each iteration, the sum of all costs ($Total\_Cost$) is updated until it converges within a defined confidence interval. At the beginning of the optimisation, all elements of $[TTD_{i,o}]$ (Time To Downtime) are set to zero, i.e., it is assumed that if there is a failure, a vehicle arrives instantly for the repair/replacement.

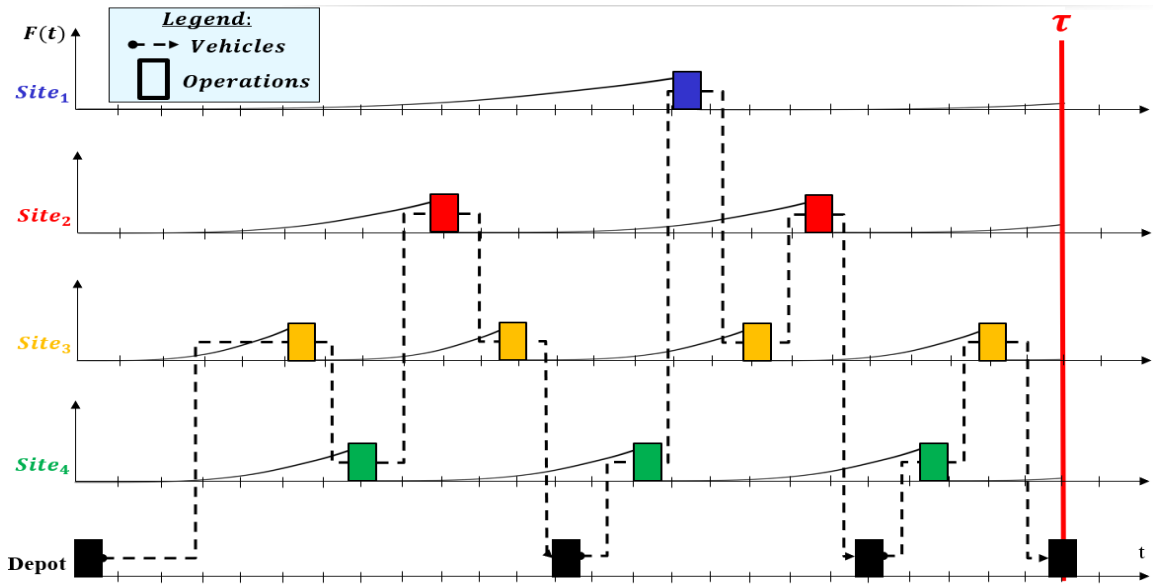

Figure 5  Example of scheduling for 4 heterogeneous sites.





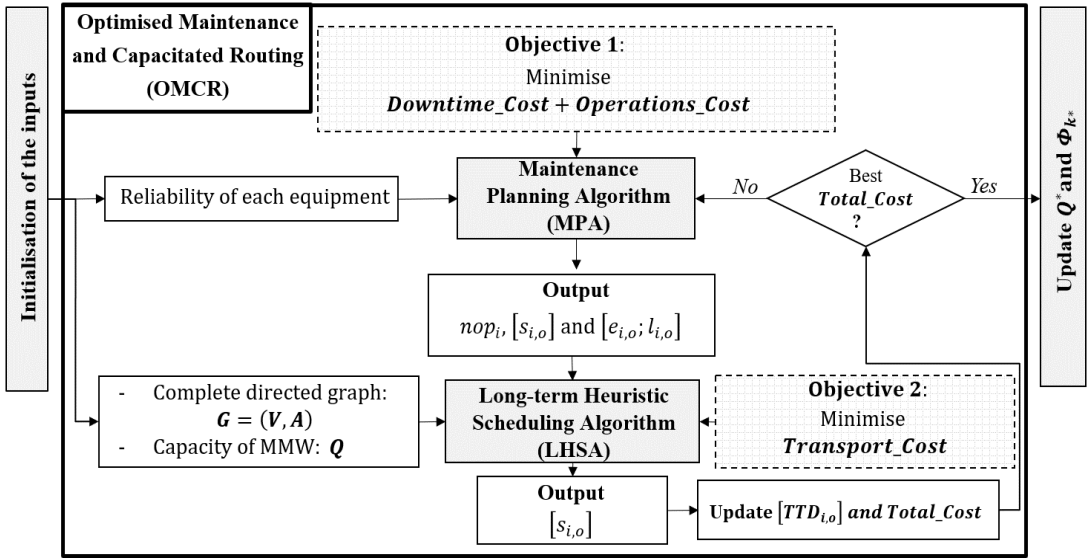

Figure 6: OMCR algorithm for combining routing and scheduling optimisation

This consideration allows the first algorithm, MPA, to propose the minimum number of maintenance operations for a site i ($nop_i = 1$). With this input data, the second algorithm, LHSA, determines the optimal route for performing the operations and, consequently, the actual value of $[TTD_{i,o}]$, which is not zero due to transportation delays and vehicle availability constraints.

$$TTD_{i,o} = (s_{i,o} - s_{i,o-1}) - \int_{s_{i,o-1}}^{s_{i,o}} \frac{t \cdot f_i(TTF_{i,o} \le t)}{F_i(TTF_{i,o} \le s_{i,o})} \quad (15)$$

Equation (15) is derived from the demonstration in Appendix B of López-Santana et al. (2016) with :

$f_i(TTF_{i,o} \le t)$ the probability density function of "Time To Failure" such that:

$$F_i(TTF_{i,o} \le s_{i,o}) = \int_{s_{i,o-1}}^{s_{i,o}} f_i(TTF_{i,o} \le t) \cdot dt \quad (16)$$

With the new $[TTD_{i,o}]$ as input, the MPA algorithm logically increases the number of operations in an attempt to reduce site unavailability, but this action contributes to increasing the cost of operations. The optimal number of operations must be found each time, followed by the optimisation of the corresponding routes. The average availability of the sites can be calculated using the simple formula (17), which is a ratio between the period during which the equipment at the site is probably in operation and the scheduling horizon.

$$A_i = \frac{\tau - \sum_{o=1:nop_i} TTD_{i,o} \cdot F_i(TTF_{i,o} \le s_{i,o})}{\tau} \quad (17)$$

After defining the OMCR algorithm used to optimise the scheduling of maintenance operations, the next step is to determine a cost-effective location for the depot.

## 4.2. Centralised Maintenance Workshop Location

Let's consider $\Omega$ as the set of all possible geographical positions for the central workshop/depot. The simplest approach to position the depot would be to evaluate the costs for all elements of $\Omega$ to compare them and deduce an optimal position. Let $\Phi\_k \in \Omega$ be the sought optimal position with k as an index that allows us to explore all elements of $\Omega$ i.e. $k \in [1; |\Omega|]$. However, the larger $\Omega$ is, the more difficult, if not impossible, it is to calculate and compare the costs of all the elements of the set. The objective is then to propose a heuristic that guarantees a small size of $\Omega$. An empirical study addresses this problem by proposing to construct $\Omega$ from sites' location (Simeu-Abazi and Gascard, 2020). The idea is to position the depot near one of the sites. However, the more sites there are, the larger the size of $\Omega$ becomes. The aim is therefore to propose a method where $|\Omega|$ is independent of the number of sites. The novelty of this paper is to position the depot at the weighted barycentre of the sites. The idea is that the depot is closer to the most critical site without being too far from the less critical ones. In this study, the criticality of each piece of equipment on site i is modelled by the cumulative distribution function ($F_i$), as presented in Section 3. The greater the probability that a piece of equipment will fail in the scheduling horizon, the closer the depot will be.

In the following, the two approaches are compared empirically through experiments. As explained in the previous paragraph, the first approach is based on the choice of the best site as a location for the depot ("*near to a site*"). The second approach is based on the barycentre of the sites whose failure probabilities are weights ("*barycenter*"). By choosing an orthonormal reference frame (0, x, y) for the geographical locations, we have:





i.   "*near to a site*": (depot(x), depot(y)) ∈ {(site$_i$ (x), ⌊site$_i$ (y))}

ii.  "*barycenter*":
$$\begin{cases} depot(x) = \sum \frac{site_i(x).F_i\,(T \le \tau)}{\sum F_i\,(T \le \tau)} \\ depot(y) = \sum \frac{site_i(y).F_i\,(T \le \tau)}{\sum F_i\,(T \le \tau)} \end{cases}$$

After researching a set of methods for cost evaluation and optimisation in a distributed maintenance context (Section 4.1) and choosing the geographical location of the central workshop/depot (Section 4.2), we will now construct a general framework to combine all the proposed methods.

### 4.3. General framework

The long-term scheduling, vehicles' capacity (Q*) and depot location $\Phi_k^*$ that minimise *Total_Cost* are obtained through the general framework illustrated in Figure 7. Indeed, it is an iterative algorithm that takes as input the geographical data of the sites ([d$_{ij}$]; [t$_{ij}$]), the health state of the equipment (F$_i$), a predefined range of vehicles' capacity [Q$_{min}$;Q$_{max}$]) and a set of probable positions of the depot (Ω).

The objective is to employ the previously developed OMCR (Optimised Maintenance and Capacitated Routing) algorithm to determine the best possible maintenance cost for each of the values of Q ∈ [Q$_{min}$ ; Q$_{max}$] and Φ_k∈Ω, through long-term scheduling. Therefore, two iterative loops are created to test all values of Q and Φ$_k$. Each iteration updates the lowest *Total_Cost* value and associated parameters

## 5. Experiments: implementation of the proposed algorithms

As a reminder, the purpose of this paper is to reduce distributed maintenance costs Long-term versus short-term scheduling: influence of the scheduling horizon

The first contribution of this paper concerns long-term scheduling as described above. The idea of this first part of the experiment is to determine the evolution of the maintenance cost (*Total_Cost*) as a function of the scheduling horizon (τ). This experiment concerns equipment whose number of failures follows the Weibull distribution with the given cumulative distribution function:

$$F_i(TTF_i \le t) = 1 - e^{-\left(\frac{t}{\eta}\right)^{\beta}} \tag{18}$$

Where η is scale parameter in years and β is the shape parameter.

It is necessary to study several scenarios for two reasons: to extend the scope of this experiment and to select the most representative scenarios for analysis. A scenario is considered relevant when it allows studying the influence of the scheduling horizon (τ), and hence long-term scheduling, i.e., when τ is large, on the performance of distributed maintenance. Assuming that, depending on the application (aerospace, railway, energy, etc.), the cost of equipment unavailability penalties can be a key variable in strategic maintenance decision-making, a scenario is represented by CP$_i$. Consequently, the next steps involve understanding the impact of CP$_i$ on the influence of τ. Thus, we vary τ from 2 months (61 days) to 2 years (732 days). The following question will guide us through this: how does the influence

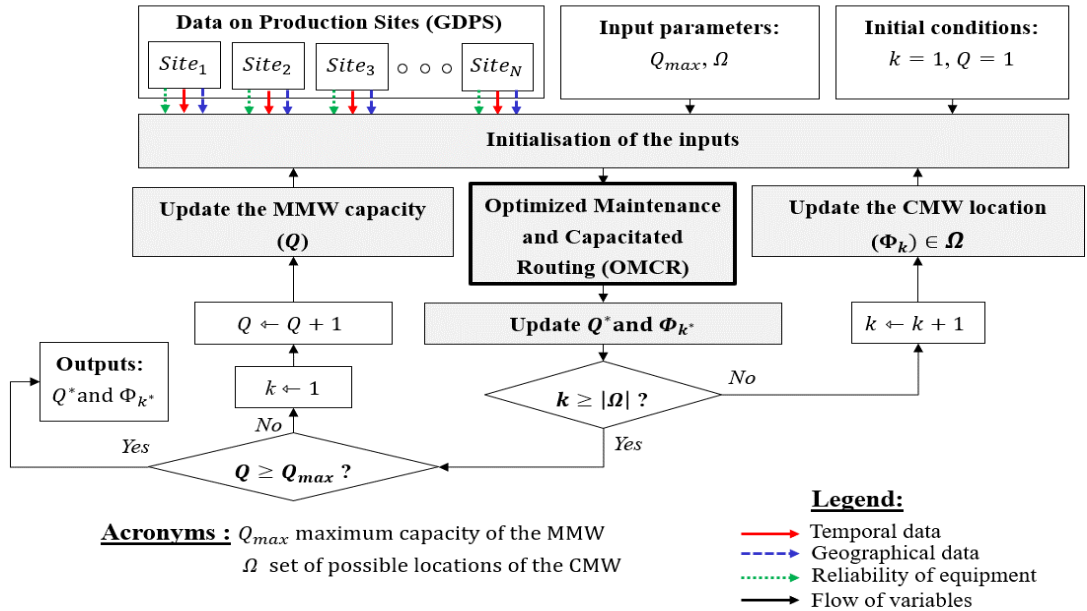

**Acronyms :** $Q_{max}$ maximum capacity of the MMW
Ω set of possible locations of the CMW

**Legend:**
— Temporal data
- - - Geographical data
··· Reliability of equipment
— Flow of variables

Figure 7: General framework of distributed maintenance optimisation.

of the scheduling horizon on the performance of distributed maintenance evolve when the cost of unavailability penalties is relatively low or high compared to the cost of replacing a failed equipment?

## 5.1. Depot location and extension of the number of sites

The second contribution of this paper concerns the location of the depot at the barycentre of the sites. The idea of this second part of the experiment is to assess whether positioning the depot at the barycentre of the sites is more interesting than the most relevant method in the literature, which positions the depot close to one of the sites. Another objective of this experiment is to study the impact of a possible extension of the number of sites after the choice of the geographical position of the depot. The idea is to evaluate the costs if additional sites were added after the construction of the depot at the proposed barycentre. The questions that will be used as a guide are as follows: is it relevant to position the depot at the barycentre of the sites? What would happen if the depot position is chosen for an initial number of sites N and then other sites are added later without changing the depot?

## 5.2. Case study - Input Data

This experiment considers 10 sites geo-distributed over a radius of 50 kilometres, corresponding approximately to an area the size of the Isère department in France.

Each site has a piece of equipment subject to uncertain failures and whose characteristics are presented in Table 3. Unlike the exponential distribution regularly used in the literature, the choice of the Weibull distribution allows us to study equipment whose failure rate varies over time. This experiment will focus on pieces of equipment that are reaching the end of their useful life and therefore require greater attention in terms of maintenance due to the increasing number of breakdowns. To consider the heterogeneity of equipment, half of the production sites have equipment that tends to deteriorate regularly ($\beta$=2) and the other half have equipment that deteriorate more rapidly ($\beta$=3). Furthermore, maintenance/replacement costs are generally high in the railway, aircraft and oil sectors. We have assumed a value of $100,000 per equipment, which is less than 0.33% of the price of a TGV (Duteil, 2016) or 0.02% of the price of an Airbus A380 (Reuters, 2019).

Only one spare part is required for each operation. The geographical positions of the production sites are randomly selected on a Cartesian plane (0,x,y) following a uniform distribution. Subsequently, the distance between each pair of equipment is calculated using the Euclidean method. For each scenario of $CP_i$ as described above, the experiment is replicated more than 10 times to ensure a 95% confidence interval for the results. Consequently, the results obtained represent the average of all replications with a 5% error rate. Appendix 1 provides an illustration of 10 replications of the positioning of production sites, noting that the positions change with each replication. A uniform distribution has been chosen to ensure that the positions of the production sites do not influence the results, as this study is solely concerned with the influence of the MMW/depot location and MMW/vehicles capacity. A fleet of homogeneous vehicles is considered and the capacity needs to be optimised. The characteristics of each vehicle are presented in Table 4. Three types of vehicles are considered, each with a nominal speed of 80 km/h.

   i.   Light: Q = 4 pieces of equipment;

   ii.   Medium: Q = 6 pieces of equipment;

   iii.   Heavy: Q = 8 pieces of equipment.

We chose Scilab 5.5.2 to implement the case study. All tests were performed using the MILP-adapted "FOSSEE Optimisation Toolbox" library. We performed the experiments on a Windows 8, 64-bit machine with an Intel(R) Core (TM) i7-10850H, 2.70 GHz CPU and 32 GB RAM. In the end, more than 1380 experiments are conducted.

Table 3: Data of the sites and equipment

| Symbols | Values | Comments | Units |
|---|---|---|---|
| $N$ | 10 | N dispersed production sites | - |
| $F_i$ $(TTF_i \leq t)$ | $1 - e^{-(\frac{t}{\eta})^\beta}$ | Weibull's cumulative distribution function with parameters $\eta$ and $\beta$ | - |
| $\eta$ | 1 | Weibull's scale parameter | $year$ |
| $\beta$ | $\beta \in \{2,3\}$ | Weibull's shape parameter | - |
| $R$ | 50 | The sites are randomly distributed in a radius $R$ | $km$ |
| $v$ | 80 | Average transport speed | $km/h$ |
| $\tau$ | $\tau \in ]0;2]$ | Horizon for maintenance operations | $year$ |
| $MTTR_i$ | 3 | Maintenance/replacement time of a piece of equipment | $h$ |
| $CR_i$ | 100,000 | Maintenance/replacement cost of a piece of equipment | $ |
| $CP_i$ | $CP_i \in [10; 100$ | Penalty cost of waiting for replacement/maintenance of a piece of equipment | $/h$ |

Table 4: data of the vehicles

| Symbols | Values | Comments | Units |
|---|---|---|---|
| $Q$ | $Q \in \{4; 6; 8\}$ | Maximum transport capacity | $Equipment$ |
| $CD$ | 2 | Cost per unit of transport distance | $/km$ |
| $CT$ | 30 | Cost per unit of transport time | $/h$ |

## 6. Results and discussion

Experimental results show that long-term scheduling (more than two months) could result in lower costs than short-term approaches (two months). With the proposed OMCR algorithm, maintenance costs can be optimised over a 2-year horizon. Therefore, different depot locations and vehicle capacities can be tested to determine the most cost-effective options. An optimal evaluation of the number of predictive maintenance operations and their scheduling could reduce total maintenance costs by up to 50%.

### 6.1. Analysis and interpretation

#### 6.1.1. Influence of penalty costs

The penalty cost represents the loss of revenue when equipment is unavailable due to failure. Figure 8 illustrates the influence of this cost on the availability of the production sites. The interpretation is similar from a cost perspective (Appendix 2). It can be seen that when $CP_i$ =\$10/h, the longer the scheduling horizon, the lower the availability (from 0.98 if $\tau$=2 months to 0.01 if $\tau$=2 years). In this case, scheduling for the long-term has no positive effect on the availability of the production sites. The penalty cost is not high enough to allow the algorithm to perform more maintenance operations if the scheduling horizon is long. In terms of cost (see Appendix 2), it is clearly observed that the larger the chosen horizon, the lower the *Total_Cost*. This demonstrates that for certain application cases where the penalty cost is not sufficiently high, corrective maintenance might be favoured, i.e., with very few preventive maintenance operations. In this scenario ($CP_i$=\$10/h), this would amount to 3 preventive operations each 2 years at a minimum *Total_Cost* of \$109/h.

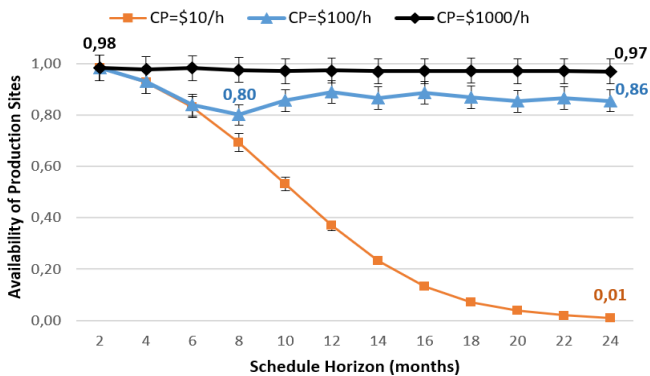

Figure 8: Influence of the penalty cost ($CP_i$)

If $CP_i$ =\$1000/h, the availability remains almost constant (0.980.05) despite the increase in the scheduling horizon. This can be explained by the fact that the penalty cost is very high. The proposed algorithm adapts the number of operations in order to maintain the production sites in an optimal availability regardless of the chosen scheduling horizon. This scenario is the one studied by (López-Santana et al., 2016). Although it guarantees a high availability of the sites, this value of the penalty cost does not highlight the influence of the choice of the scheduling horizon. In terms of cost (see Appendix 2), we naturally observe a variation of less than 6% in *Total_Cost*, regardless of whether the scheduling horizon is short or long. It might indicate that for an application domain where $CP_i$ =\$1000/h (high penalty cost in case of failure), long-term scheduling would not have a very significant impact.

If $CP_i$ =\$100/h then two phases of availability evolution are observed. From 2 months to 8 months horizon, the availability decreases from 0.98 to 0.80 and then slowly increases from 0.80 to 0.86. In this case, the value of the penalty cost is interesting for the observation of the scheduling horizon influence. The following sections will provide a more in-depth interpretation of this scenario.

As a reminder, a scenario can be considered specific to a particular application domain. For example, in the aerospace domain, the cost of unavailability is very high (in this case, setting $CP_i$ =\$1000/h would be appropriate). In the railway domain, the cost of unavailability is generally lower than in the aerospace domain (setting $CP_i$ =\$100/h would be interesting in this case). In other industrial domains where equipment unavailability does not lead to significant costs, setting $CP_i$ =\$10/h would be acceptable. We will now focus on the scenario with $CP_i$ =\$100/h, as it offers more insight into long-term scheduling, as observed earlier.

#### 6.1.2. Optimal scheduling horizon

Figure 9 shows the evolution of costs according to the scheduling horizon for the scenario $CP_i$ =\$100/h. The 4 types of costs can be distinguished: *Transport_Cost*, *Operations_Cost*, *Downtime_Cost* and *Total_Cost* Transport costs are in the minority compared to other costs.

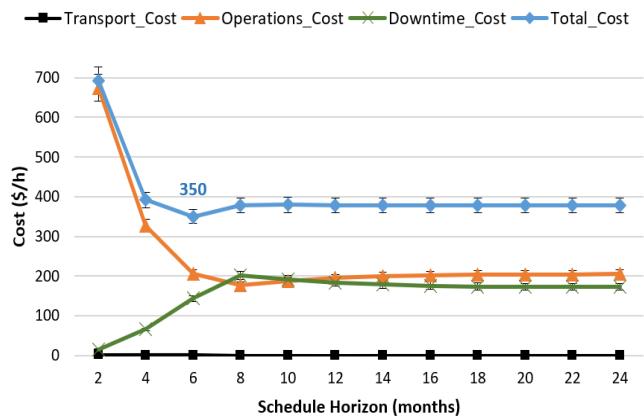

Figure 9: Costs and schedule horizon

The two phases discovered above (decrease phase and slow growth phase) have different effects on costs. From 2 to 8 months the costs of operations decrease while the costs of downtime increase, then the trends are reversed. These evolutions allow the total costs to have a local minimum for a 6-month scheduling horizon (from $700/h to $350/h), a 50% decrease. Therefore, for this scenario it is effective from a cost point of view to plan the maintenance operations every 6 months. After determining an optimal time horizon for maintenance scheduling, the next step is to determine the location of the depot and the capacity of the transport vehicles.

### 6.1.3. Influence of depot location

Figure 10 shows the influence of the geographical position of the depot on the distance travelled by the vehicles (a) and the calculation time of the evaluation algorithm (b). The location of the depot was determined using the proposed general framework, for 10 production sites. Several other sites (up to 30 sites) were then added without changing the location of the depot. The first obvious observation is that the more production sites there are, the further the vehicles travel each year. Secondly, as the number of sites increases, the positioning of the depot at the barycentre of the sites becomes more advantageous. These results show that it could be more interesting to position the depot at the barycentre than to position it near one of the production sites. This choice is even more advantageous if the number of sites increases in time after the installation of the depot.

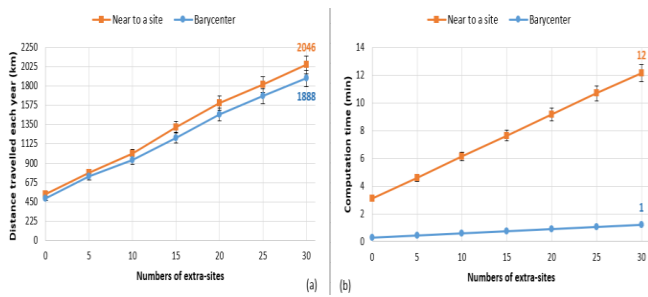

Figure 10: Influence of the depot location

### 6.1.4. Influen5ce of vehicle capacity

Figure 11 shows the influence of vehicle capacity on annual distance travelled (a) and transport costs (b). As in the interpretation in the previous section, the more additional production sites there are, the greater the annual distance travelled. It can be seen that heavy vehicles travel less distance than light vehicles (1499km and 1973km), a difference of 474km. In fact, light vehicles carry less equipment than heavy vehicles and therefore have to travel more kilometres to meet the demand. But as far as costs are concerned, it is rather the light vehicles that are less expensive. Indeed, each kilometre driven by a heavy vehicle is more expensive than that of a light vehicle. Therefore, light vehicles could be chosen to ensure the maintenance of the sites at a lower cost. Of course, these conclusions have been drawn for this particular case study. The main contributions remain the model and the optimisation methodology.

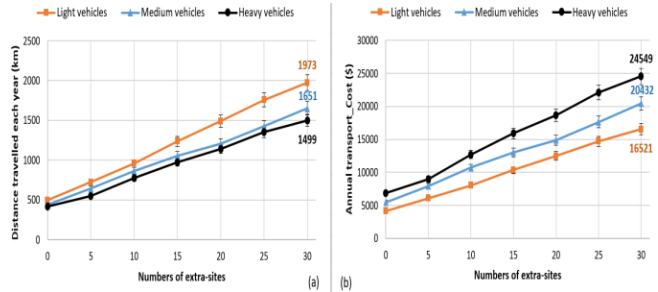

Figure 11: Influence of vehicle capacity

## 6.2. Industrial usefulness

As presented in the introduction, distributed maintenance can be applied in several areas: oil & gas, railway and aircraft domains. The key is to have several geographically distributed production sites and a centralised entity responsible for the maintenance of all sites. The results of the experiments carried out highlight the relevance of the proposed approach for the implementation of distributed maintenance. The case of 10 production sites has been studied and our approach can be applied to many more (30 additional sites were tested) thanks to the heuristics proposed to solve the computation time problem.

The production sites were considered heterogeneous, i.e. each site has different equipment. This assumption is close to the industrial reality where the geo-distributed sites generally have different facilities or different states of health. However, this study assumes that each site has only one piece of production equipment. In other applications where there is more than one piece of equipment per site, this method could be adapted by considering only the most critical equipment per site.

The vehicle fleet is considered homogeneous, i.e. all vehicles have the same transport capacity. Although this assumption allows the capacity constraint of the mobile maintenance workshop fleet to be considered, it remains limited for certain applications. The proposed method could then be applied separately to each type of vehicle. A weighted mean of the costs could then be used to obtain an approximate result in the absence of a more advanced method.

In this study, the results show that it is possible to estimate maintenance costs over a long-time horizon. This approach would allow practitioners to compare different maintenance policies over the long term. A broad view in time offers the possibility to make effective decisions long before failures occur.

## 6.3. Conclusion

This study aimed to design and optimize distributed maintenance operations for a network of geographically dispersed production sites requiring regular equipment servicing. The proposed strategy involved centralizing maintenance activities within a dedicated workshop and deploying a fleet of vehicles to perform on-site interventions. The originality of the work lies in its focus on three key decision-making areas: the long-term scheduling of maintenance operations, the optimal geographical positioning of the central maintenance facility, and the strategic sizing of the vehicle fleet.

To address these challenges, a hybrid approach combining a linear optimization model with heuristic algorithms was developed to minimize overall maintenance costs. A case study demonstrated that systematic long-term scheduling can reduce maintenance expenditures by up to 50%, primarily through precise estimation of required maintenance actions to maintain site availability over an extended planning horizon. The equipment considered followed a Weibull failure distribution with a one-year scale parameter. In this context, scheduling interventions every two months proved cost-effective when penalty costs for downtime were high relative to replacement costs; otherwise, a six-month interval was more optimal.

Locating the central maintenance workshop at the weighted barycentre of the production sites was shown to significantly reduce computational complexity—by a factor proportional to the number of sites involved. While light vehicles covered longer distances than heavy vehicles, their lower operational costs made configurations involving multiple light vehicles (e.g., two vehicles carrying four units each) more economically viable than fewer, higher-capacity heavy vehicles.

Nevertheless, the study has several limitations, both in terms of modelling and optimization methodology. First, it assumes that each repair restores the equipment to an "as good as new" state, overlooking scenarios involving imperfect repairs. Future research could investigate how varying levels of repair quality influence long-term maintenance costs. Second, the cost function was decomposed into downtime and transportation components, each addressed through separate optimization routines. A promising direction would be to develop a unified model that jointly optimizes both components, for instance, by integrating them into a single objective function via a linear formulation. Third, we considered only one piece of equipment per production site. Generalization to multiple pieces of equipment with multiple failures is currently under study, but the model and its resolution remain complex using linear programming.

Further extensions could include the incorporation of unexpected equipment failures into the scheduling model, along with the development of online algorithms capable of dynamically adjusting vehicle routes in response. Moreover, the current study focuses solely on cost as the performance metric. Future investigations could broaden the scope by integrating additional criteria, such as $CO_2$ emissions from transportation activities. Another potential research avenue involves the dimensioning of spare parts inventory at the central workshop, which was assumed to have unlimited capacity in this work. Finally, applying the proposed methodology to a real industrial case study would enhance the practical relevance of the findings and help address data availability challenges encountered in this initial exploration. From a broader perspective, this hybrid approach addresses a key limitation in traditional PHM systems by providing both generalization and physical interpretability — two often competing objectives. Moreover, the framework is modular and scalable, making it adaptable to a wide range of industrial domains, including aerospace, energy systems, and complex manufacturing processes.

Several avenues of research are identified and need to be addressed in the immediate future. First, real-time deployment of the hybrid model will require efficient online learning techniques and adaptive filtering to process streaming data. Second, further research will focus on generalizing the framework to multi-failure scenarios and system-wide degradation mechanisms, using probabilistic graphical models or physics-based neural networks. Finally, integrating Explainable AI (XAI) with ML components would improve transparency and user trust, facilitating wider adoption of hybrid PHM systems in safety-critical environments.

### ACKNOWLEDGEMENT

This work is supported by the National Research Agency as part of Investments for the Future (ANR-15-IDEX-02) - CIRCULAR cross-cutting program.

### REFERENCES

Allaham, H., Dalalah, D., 2022. MILP of multitask scheduling of geographically distributed maintenance tasks. Int. J. Ind. Eng. Comput. 13, 119–134.

Borcinova, Z., 2017. Two models of the capacitated vehicle routing problem. Croat. Oper. Res. Rev. 8, 463–469. https://doi.org/10.17535/crorr.2017.0029

Cakirgil, S., Yucel, E., Kuyzu, G., 2020. An integrated solution approach for multi-objective, multi-skill workforce scheduling and routing problems. Comput. Oper. Res. https://doi.org/10.1016/j.cor.2020.104908

Campbell, S. L., Chancelier, J. P., Nikoukhah, R., (2010). Modeling and Simulation in SCILAB (pp. 73-106). Springer New York.

Camci, F., 2015. Maintenance scheduling of geographically distributed assets with prognostics information. Eur. J. Oper. Res. https://doi.org/10.1016/j.ejor.2015.03.023

Djeunang Mezafack, R.A., Di Mascolo, M., Simeu-Abazi, Z., 2022. Systematic literature review of repair shops: focus on

sustainability. Int. J. Prod. Res. 60, 7093–7112. https://doi.org/10.1080/00207543.2021.2002965

Djeunang Mezafack R. A., Simeu-Abazi Z., Di Mascolo M., "Optimisation of Distributed Maintenance: Design, Scheduling and Capacitated Routing Problem". Advanced Maintenance Engineering, Services and Technologies - 5th AMEST 2022 Bogota

Djeunang Mezafack R.A. La maintenance distribuée au cœur d'une économie circulaire : contribution à sa mise en œuvre et à son évaluation. Automatique / Robotique. PhD Université Grenoble Alpes 2023 HAL Id:tel-04514328.

Drake, J.H., Starkey, A., Owusu, G., Burke, E.K., 2020. Multiobjective evolutionary algorithms for strategic deployment of resources in operational units. Eur. J. Oper. Res. 282, 729–740. https://doi.org/10.1016/j.ejor.2019.02.002

Duteil, M., 2016. 15 TGV pour Alstom : une commande au prix fort pour l'État [WWW Document], URL https://www.europe1.fr/economie/15-tgv-pour-alstom-une-commande-au-prix-fort-pour-letat-2864806 (accessed 6.19.23).

Fontecha, J.E., Guaje, O.O., Duque, D., Akhavan-Tabatabaei, R., Rodriguez, J.P., Medaglia, A.L., 2020. Combined maintenance and routing optimization for large-scale sewage cleaning. Ann. Oper. Res. https://doi.org/10.1007/s10479-019-03342-8

Gopalakrishnan, M., Subramaniyan, M., Skoogh, A., 2022. Data-driven machine criticality assessment – maintenance decision support for increased productivity. Prod. Plan. Control 33, 1–19. https://doi.org/10.1080/09537287.2020.1817601

Gupta, D., 2003. A new algorithm to solve Vehicle Routing Problems (VRPs). Int. J. Comput. Math. https://doi.org/10.1080/0020716022000005537

Hani, Y., Amodeo, L., Yalaoui, F., Chen, H., 2007. Ant colony optimization for solving an industrial layout problem. Eur. J. Oper. Res. 183, 633–642. https://doi.org/10.1016/j.ejor.2006.10.032

Hedjazi, D., Layachi, F., Boubiche, D.E., 2019. A multi-agent system for distributed maintenance scheduling. Comput. Electr. Eng. 77, 1–11. https://doi.org/10.1016/j.compeleceng.2019.04.016

Jia, C., Zhang, C., 2020. Joint optimization of maintenance planning and workforce routing for a geographically distributed networked infrastructure. IISE Trans. https://doi.org/10.1080/24725854.2019.1647478

Konstantakopoulos, G.D., Gayialis, S.P., Kechagias, E.P., 2020. Vehicle routing problem and related algorithms for logistics distribution: a literature review and classification. Oper. Res. https://doi.org/10.1007/s12351-020-00600-7

López-Santana, E., Akhavan-Tabatabaei, R., Dieulle, L., Labadie, N., Medaglia, A.L., 2016. On the combined maintenance and routing optimization problem. Reliab. Eng. Syst. Saf. 145, 199–214. https://doi.org/10.1016/j.ress.2015.09.016

Manco, P., Rinaldi, M., Caterino, M., Fera, M., Macchiaroli, R., 2022. Maintenance management for geographically distributed assets: a criticality-based approach. Reliab. Eng. Syst. Saf. 218, 108148. https://doi.org/10.1016/j.ress.2021.108148

Mariescu-Istodor, R., Cristian, A., Negrea, M., Cao, P., 2021. VRPDiv: A Divide and Conquer Framework for Large Vehicle Routing Problems. ACM Trans. Spat. Algorithms Syst. 7, 23:1-23:41. https://doi.org/10.1145/3474832

Meng, F.C., 2000. Relationships of Fussell–Vesely and Birnbaum importance to structural importance in coherent systems. Reliab. Eng. Syst. Saf. 67, 55–60. https://doi.org/10.1016/S0951-8320(99)00043-5

Ndiaye I. 2014. Maintenance distribuee : mise en oeuvre d'un algorithme d'optimisation des coûts de maintenance Master report, June 2014.

Nguyen, H.S.H., Do Van, P., Vu, H.C., Iung, B., 2019. Dynamic maintenance grouping and routing for geographically dispersed production systems. Reliab. Eng. Syst. Saf. 185, 392–404. https://doi.org/10.1016/j.ress.2018.12.031

Rashidnejad, M., Ebrahimnejad, S., Safari, J., 2018. A bi-objective model of preventive maintenance planning in distributed systems considering vehicle routing problem. Comput. Ind. Eng. 120. https://doi.org/10.1016/j.cie.2018.05.001

Razavi Al-e-hashem, S.A., Papi, A., Pishvaee, M.S., Rasouli, M., 2022. Robust maintenance planning and scheduling for multi-factory production networks considering disruption cost: a bi-objective optimization model and a metaheuristic solution method. Oper. Res. 22, 4999–5034. https://doi.org/10.1007/s12351-022-00733-x

Reuters, 2019. Airbus devra baisser le prix de l'A380 pour en vendre plus [WWW Document], Les Echos Investir. URL https://investir.lesechos.fr/actu-des-valeurs/la-vie-des-actions/airbus-devra-baisser-le-prix-de-la380-pour-en-vendre-plus-iag-1825402 (accessed 6.19.23).

Saihi, A., Ben-Daya, M., As'ad, R.A., 2022. Maintenance and sustainability: a systematic review of modeling-based literature. J. Qual. Maint. Eng. ahead-of-print. https://doi.org/10.1108/JQME-07-2021-0058

Sanchez, D.T., Boyacı, B., Zografos, K.G., 2020. An optimisation framework for airline fleet maintenance scheduling with tail assignment considerations. Transp. Res. Part B Methodol. 133, 142–164. https://doi.org/10.1016/j.trb.2019.12.008

Sedghi, M., Kauppila, O., Bergquist, B., Vanhatalo, E., Kulahci, M., 2021. A taxonomy of railway track maintenance planning and scheduling: A review and research trends. Reliab. Eng. Syst. Saf. 215, 107827. https://doi.org/10.1016/j.ress.2021.107827

Si, G., Xia, T., Pan, E., Xi, L., 2022. Service-oriented global optimization integrating maintenance grouping and technician routing for multi-location multi-unit production systems. IISE Trans. 54, 894–907. https://doi.org/10.1080/24725854.2021.1957181

Simeu-Abazi, Z., Ahmad, A.A., 2011. Optimisation of distributed maintenance: Modelling and application to the multi-factory production. Reliab. Eng. Syst. Saf. 96, 1564–1575. https://doi.org/10.1016/j.ress.2011.05.011

Simeu-Abazi, Z., Gascard, E., 2020. Implementation of a cost optimization algorithm in a context of distributed maintenance, in: 2020 International Conference on Control, Automation and Diagnosis (ICCAD). Presented at the 2020 International Conference on Control, Automation and Diagnosis (ICCAD), pp. 1–6. https://doi.org/10.1109/ICCAD49821.2020.9260507

Sleptchenko, A., Turan, H.H., Pokharel, S., ElMekkawy, T.Y., 2019. Cross-training policies for repair shops with spare part inventories. Int. J. Prod. Econ. 209, 334–345. https://doi.org/10.1016/j.ijpe.2017.12.018

Tang, H., Miller-Hooks, E., Tomastik, R., 2007. Scheduling technicians for planned maintenance of geographically distributed equipment. Transp. Res. Part E Logist. Transp. Rev. 43, 591–609. https://doi.org/10.1016/j.tre.2006.03.004

Valet, A., Altenmüller, T., Waschneck, B., May, M.C., Kuhnle, A., Lanza, G., 2022. Opportunistic maintenance scheduling with deep reinforcement learning. J. Manuf. Syst. 64, 518–534. https://doi.org/10.1016/j.jmsy.2022.07.016

Wang, K., Djurdjanovic, D., 2018. Joint Optimization of Preventive Maintenance, Spare Parts Inventory and Transportation Options for Systems of Geographically Distributed Assets. MACHINES. https://doi.org/10.3390/machines6040055

Wu, S., Coolen, F.P.A., 2013. A cost-based importance measure for system components: An extension of the Birnbaum importance. Eur. J. Oper. Res. 225, 189–195. https://doi.org/10.1016/j.ejor.2012.09.034

Yulong, L., Chi, Z., Chuanzhou, J., Xiaodong, L., Yimin, Z., 2019. Joint optimization of workforce scheduling and routing for restoring a disrupted critical infrastructure. Reliab. Eng. Syst. Saf. https://doi.org/10.1016/j.ress.2019.106551

Zhang, C., Yang, T., 2021. Optimal maintenance planning and resource allocation for wind farms based on non-dominated sorting genetic algorithm-II. Renew. Energy 164, 1540–1549. https://doi.org/10.1016/j.renene.2020.10.125.